\documentclass[authoryear,11pt,preprint]{elsarticle}

\usepackage{color}
\usepackage{lineno}
\usepackage{pdflscape}
\usepackage{amsmath}
\usepackage{booktabs}
\usepackage[hyphens]{url}
\usepackage{hyperref}
\usepackage{multirow}
\usepackage{addlines}
\usepackage{adjustbox}
\usepackage{natbib} 
\usepackage[english]{babel}
\usepackage[dvipsnames]{xcolor}
\usepackage{soul}
\usepackage{adjustbox}
\usepackage{verbatim}
\usepackage[margin=2.5cm]{geometry}
\usepackage{mathtools}
\usepackage{tikz}
\usepackage{pgfplots}
\usepackage{enumitem}
\usepackage{calc}

\definecolor{darkblue}{rgb}{0,0.2,0.4}
\usepackage{graphicx}
\usepackage{subcaption}
\usepackage{amssymb}
\usepackage{amsthm}
\usepackage{etoolbox}
\usepackage{algorithm}
\usepackage{algpseudocode}

\makeatletter
\newcommand*{\algrule}[1][\algorithmicindent]{%
  \hspace*{.2em}
  \vrule 
  \hspace*{\dimexpr#1-.2em-.4pt}%
}

\newcommand{\StatePar}[1]{%
  \State\parbox[t]{\dimexpr\linewidth-\ALG@thistlm}{\strut #1\strut}%
}
\renewcommand{\ALG@beginalgorithmic}{\offinterlineskip}

\newcount\ALG@printindent@tempcnta
\def\ALG@printindent{%
  \ifnum \theALG@nested > 0
    \ifx\ALG@text\ALG@x@notext
    \else
      \unskip
      \ALG@printindent@tempcnta=1
      \loop
        \algrule[\csname ALG@ind@\the\ALG@printindent@tempcnta\endcsname]%
        \advance \ALG@printindent@tempcnta 1
        \ifnum \ALG@printindent@tempcnta<\numexpr\theALG@nested+1\relax
      \repeat
        \fi
    \fi
}
\patchcmd{\ALG@doentity}{\noindent\hskip\ALG@tlm}{\ALG@printindent}{}{\errmessage{failed to patch}}
\makeatother

\algrenewcommand\algorithmicend{\strut\textbf{end}}
\algrenewcommand\algorithmicdo{\strut\textbf{do}}
\algrenewcommand\algorithmicwhile{\strut\textbf{while}}
\algrenewcommand\algorithmicfor{\strut\textbf{for}}
\algrenewcommand\algorithmicforall{\strut\textbf{for all}}
\algrenewcommand\algorithmicloop{\strut\textbf{loop}}
\algrenewcommand\algorithmicrepeat{\strut\textbf{repeat}}
\algrenewcommand\algorithmicuntil{\strut\textbf{until}}
\algrenewcommand\algorithmicprocedure{\strut\textbf{procedure}}
\algrenewcommand\algorithmicfunction{\strut\textbf{function}}
\algrenewcommand\algorithmicif{\strut\textbf{if}}
\algrenewcommand\algorithmicthen{\strut\textbf{then}}
\algrenewcommand\algorithmicelse{\strut\textbf{else}}
\algrenewcommand\algorithmicrequire{\strut\textbf{Input:}}
\algrenewcommand\algorithmicensure{\strut\textbf{Output:}}
\let\oldState\State
\renewcommand{\State}{\oldState\strut}
\algnewcommand{\IIf}[1]{\State\algorithmicif\ #1\ \algorithmicthen}
\algnewcommand{\EndIIf}{\unskip\ \algorithmicend\ \algorithmicif}

\usepackage{nomencl}
\makenomenclature
\usepackage{framed} 
\usepackage{multicol} 

\usepackage[utf8]{inputenc}
\usepackage{newunicodechar}

\DeclareUnicodeCharacter{0229}{\c{e}}
\DeclareUnicodeCharacter{2009}{ °}
\usepackage[T1]{fontenc}
\usepackage{textcomp}
\usepackage{lmodern}
\usepackage{ltxcmds}
\usepackage{gensymb}

\usepackage{siunitx}
\DeclareSIUnit{\euro}{\texteuro}
\makeatother
\definecolor{royalpurple}{rgb}{0.58, 0.44, 0.86}

\newcommand{\blue}[1]{{\color{blue} #1}}

\allowdisplaybreaks

\usepackage{setspace}
\hypersetup{colorlinks = true,linkcolor = blue,anchorcolor =red,citecolor = blue,filecolor = red,urlcolor = red, pdfauthor=author}

\journal{EURO Journal on Computational optimisation}

\pgfplotsset{compat=1.17} 
\begin{document}

\begin{frontmatter}

\title{Neural networks for multi-horizon stochastic programming}


\author[a]{Hongyu Zhang\corref{cor1}}
\ead{hongyu.zhang@soton.ac.uk}
\author[b]{Gabriele Sormani}
\ead{g.sormani6@campus.unimib.it}
\author[b]{Enza Messina}
\ead{enza.messina@unimib.it}
\author[c]{Alan King}
\ead{alan.king@ars-proba.ca}
\author[d]{Francesca Maggioni}
\ead{francesca.maggioni@unibg.it}

\cortext[cor1]{Corresponding author}

\address[a]{School of Mathematical Sciences, University of Southampton, Building 54, Highfield Campus, Southampton, SO14 3ZH, United Kingdom}
\address[b]{Department of Informatics, Systems and Communication, University of Milano-Bicocca, Viale Sarca 336, Milano, 20126, Italy}
\address[c]{Ars Probabilitas}
\address[d]{Department of Management, Information and Production Engineering, University of Bergamo, Via Marconi - 5 24044 Dalmine, Italy}

\begin{abstract}

This paper proposes a machine-learning-based solution approach for solving multi-horizon stochastic programs. The approach embeds a deep learning neural network into a multi-horizon stochastic program to approximate the recourse operational objective function. The proposed approach is demonstrated on a UK power system planning problem with uncertainty at investment and operational timescales. The results show that (1) the surrogate neural network performs well across three different architectures, (2) the proposed approach is up to 34.72 times faster than the direct solution of the monolithic deterministic equivalent counterpart, (3) the surrogate-based solutions yield comparable in-sample stability and improved out-of-sample performance relative to the deterministic equivalent, indicating better generalisation to unseen scenarios. The main contributions of the paper are: (1) we propose a machine-learning-based framework for solving multi-horizon stochastic programs, (2) we introduce a neural network embedding formulation tailored to multi-horizon stochastic programs with continuous first-stage decisions and fixed scenario sets, extending existing surrogate modelling approaches from two-stage to multi-horizon settings, and (3) we provide an extensive computational study on a realistic UK power system planning problem, demonstrating the trade-off between approximation accuracy, computational efficiency, and solution robustness for different neural network architectures and scenario set sizes.

\end{abstract}

\end{frontmatter}


\section{Introduction}
\label{sec:introduction}
Long-term investment planning problems often face uncertainty at investment and operational timescales. Including multi-timescale uncertainty in a multistage stochastic program can easily lead to an intractable model. Multi-Horizon Stochastic Programming (MHSP) has been proposed to model multi-timescale uncertainty efficiently by reducing the model size at a cost of losing some information \citep{Kaut2014}. In essence, MHSP is a specific case of multistage stochastic programs with block separable recourse \cite{Louveaux1986MultistageRecourse}. Despite its advantage in reducing model size, some monolithic MHSP models have been intractable \citep{Zhang2025IntegratedDecomposition, Zhang2024AUncertainty} or required a lot of time to solve \citep{Backe2022EMPIRE:Analyses}. Therefore, addressing the computational difficulties in solving MHSP models has been studied \citep{Zhang2024DecompositionProgramming,Zhang2024AUncertainty,Zhang2025IntegratedDecomposition,zhang2025multi-timescale, mazzi2025adaptivebendersdecompositionenhanced}. The existing literature has mainly explored decomposing MHSP and required MHSP to have a certain structure. However, leveraging machine learning in solving MHSP efficiently has not been studied. 

In this paper, we first propose a machine-learning-based solution approach for efficiently solving MHSP problems. Machine learning techniques, and in particular neural networks, have recently shown strong potential in approximating complex mappings and accelerating the solution of high-dimensional optimisation problems \citep{patel2022neur2sp,Chou1}. In this context, we employ surrogate models, i.e., data-driven approximations of computationally expensive subproblems, to enhance the tractability of MHSP formulations. The proposed framework aims to find solutions that remain flexible and robust across possible future investment and operational scenarios. Specifically, we analyse the ability of neural network surrogates to generalise to unseen or unexpected scenarios without overfitting to the limited set of scenarios used during model training. The algorithm consists of three main steps: (1) training a surrogate neural network to approximate the operational subproblems; (2) embedding the trained neural network within the original MHSP problem; and (3) solving the resulting neural network–embedded model. The embedding procedure extends the approach proposed by \cite{patel2022neur2sp} for two-stage stochastic programming to the MHSP case. We apply the proposed methodology to a UK power system investment planning problem under long-term and short-term uncertainty. We then compare the surrogate-based solutions with the deterministic equivalent model in terms of solution quality, stability, and in-sample and out-of-sample performance.

The main contributions of this paper are summarised as follows:
\begin{itemize}
\item We first propose a machine-learning-based framework for solving MHSP problems, leveraging surrogate neural networks to approximate computationally expensive operational subproblems.
\item  We first introduce a neural network embedding formulation for MHSP, extending surrogate modelling approaches to handle multi-timescale uncertainty and ensuring scalability to large-scale instances.
\item We conduct extensive computational experiments on a UK power system investment planning problem, comparing the surrogate-based approach with the deterministic equivalent model in terms of computational efficiency, solution quality, and in-sample and out-of-sample performance.
\end{itemize}

The outline of the paper is as follows: Section \ref{sec:literature_review} introduces the literature review on solution methods for MHSP, applications of MHSP in energy system planning, and surrogate models for stochastic programming. Section \ref{sec:MHSP formulation} presents a mathematical formulation of MHSP. Section \ref{sec:surrogateNN_MHSP} mathematically introduces the neural network embedding in MHSP. Section \ref{sec:problem_description_modelling_assumptions} provides the problem description, scenario generation and assumptions. Section \ref{sec:model} presents the mathematical model for the UK power system investment planning. Section \ref{sec:results} reports the computational results and numerical analysis. Section \ref{sec:discussion} discusses the implications of the method and results and summarises the limitations of the research. Section \ref{sec:conclusions} concludes the paper and suggests further research.

\section{Literature review}
\label{sec:literature_review}
In the following, we present a brief overview of relevant literature on MHSP and corresponding solution methods, the application of MHSP in the energy sector, and surrogate models in stochastic programming problems. 

\subsection{MHSP}
\label{sec:multi_horizon_programming}
\begin{figure}[htb!]
    \centering
    \includegraphics[scale=0.9]{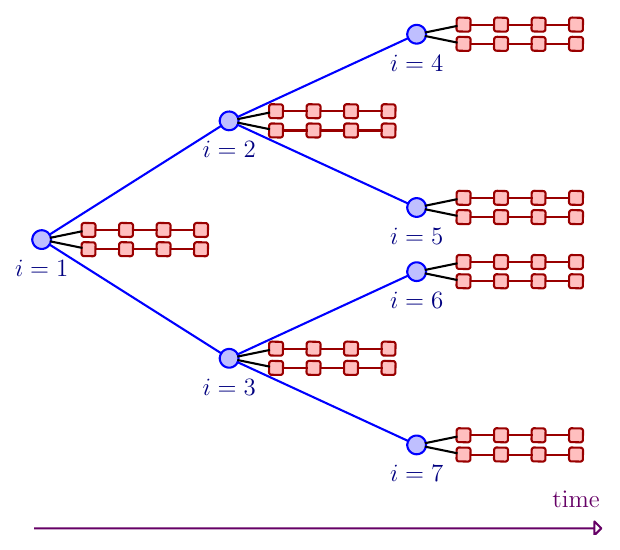}
    \caption{Illustration of MHSP with short-term and long-term uncertainty (blue circles: investment nodes, red squares: operational periods, $i:$ index of the investment nodes).}
    \label{fig:multi-horizon_tree}
\end{figure}
Investment planning of an energy system often faces uncertainty from two time horizons \citep{Kaut2014,lara2020,micheli2025}: (a) the uncertainty from the operational time horizon, such as the availability of renewable energies. The operational uncertainty becomes even more crucial for a system with higher penetration of intermittent renewable energies, and (b) the uncertainty from the strategic time horizon, e.g., long-term demand. One can have a much smaller model by disconnecting operational nodes between successive planning stages and embedding them into their respective strategic nodes. We call an operational problem embedded in a strategic node an operational node. The resulting model is called MHSP. An illustration of MHSP with short-term and long-term uncertainty is shown in Figure \ref{fig:multi-horizon_tree}. MHSP can be understood as an approximation to multistage stochastic programming unless two requirements are met \citep{Kaut2014}: (a) strategic and operational uncertainties are independent, and the strategic decisions must not depend on any particular operational decisions, and (b) the operational decisions in the last operational period in a stage do not affect the system operation in the first operational period in the next stage. The bounds in MHSP have been studied in \citep{Maggioni2020BoundsPrograms}. In this paper, we propose a surrogate neural network to approximate the objectives of the operational problems.

Despite MHSP leading to a smaller model compared with multistage stochastic programs, the computational difficulty needs to be addressed for large problems. The special structure of MHSP allows decomposing a problem with Benders-type algorithms. \cite{Zhang2024DecompositionProgramming} first established different decomposition algorithms for MHSP. A stabilised adaptive Benders decomposition algorithm was proposed in \cite{Zhang2024AUncertainty} and demonstrated on a power system investment planning problem with up to 1 billion variables and 4.5 billion constraints. The algorithm showed a significant reduction in computational time. \cite{Zhang2025IntegratedDecomposition} extended the approach in \cite{Zhang2024AUncertainty} and solved a problem with binary variables in the investment planning part in order to capture the economic scale and model retrofit and abandonment decisions. \cite{mazzi2025adaptivebendersdecompositionenhanced} proposed Adaptive Benders decomposition and enhanced SDDP for problems with convex multistage stochastic subproblems with a common structure, which can solve MHSP with multistage recourse. The proposed algorithms are efficient but they usually require the subproblems to have a certain structure. However, the proposed neural network approach in this paper does not require specific subproblem structure and still provides an efficient way to solve MHSP.

\subsection{Energy system planning}
MHSP has been applied to a variety of energy system planning problems \citep{Skar2016,Backe2022EMPIRE:Analyses,Zhang2022OffshoreShelf,Durakovic2023PoweringPrices,zhang2025multi-timescale}. Capacity expansion planning typically considers either an existing system with historical assets and operational patterns or a completely new system, and determines investment decisions required to meet demand while satisfying, among other factors, environmental constraints. Such problems are formulated using either deterministic optimisation models \citep{Lara2018DeterministicAlgorithm} or stochastic programming approaches \citep{Backe2022EMPIRE:Analyses, Conejo2016}. Although \citet{Backe2022EMPIRE:Analyses} adopted an MHSP formulation, long-term uncertainty was not incorporated. Mixed-Integer Linear Programming (MILP) is also frequently employed \citep{lara2020} to represent discrete investment decisions. Obtaining meaningful environmental and economic insights often requires solving large-scale models, as demonstrated in \citet{Li2022} and \citet{zhang2022_OEH}. To address computational challenges, \citet{Munoz2016} introduced a novel bounding scheme combined with Benders decomposition to solve a large MILP-based investment planning problem. Beyond single-carrier systems such as power, gas or heat, integrated multi-carrier energy systems have also been studied. For instance, energy hubs that convert, process, and store multiple energy carriers within an investment planning framework were examined by \citet{Zhang2024AUncertainty}. Stochastic programming has similarly been applied to natural gas systems \citep{Fodstad2016}, offshore oil and gas infrastructure planning \citep{Gupta2014}, and hydrogen networks \citep{Galan2019}. For a comprehensive review of capacity expansion planning, we refer the reader to \citet{Krishnan2016Co-optimizationApproaches}.

\subsection{Surrogate models for stochastic programming}
In this paper, we embed a deep learning model as a surrogate function to approximate the operational subproblem costs. In the following, we present a brief literature review on deep learning models and the use of surrogate models in stochastic programming. 

Deep learning models have demonstrated remarkable success in a variety of domains, including computer vision, natural language processing, and speech recognition, owing to their capacity to model complex, high-dimensional data distributions. Key architectures include Convolutional Neural Networks (CNNs), Recurrent Neural Networks (RNNs), and transformer-based models, each tailored to exploit specific structural properties of the data. 

In recent years, the use of surrogate models has emerged as an effective approach to accelerate the solution of optimisation problems. Surrogate models replace computationally expensive optimisation models with simpler, data-driven approximations that preserve key structural characteristics while significantly reducing computational cost. These models, often constructed using techniques such as polynomial regression, Gaussian processes, or neural networks, enable rapid evaluation of candidate solutions and efficient exploration of the solution space.

\cite{Fischetti2018} pioneered the explicit embedding of trained neural networks into MILP formulations, providing a framework to integrate predictive models directly within optimisation problems. Later efforts leveraged this formulation to combine deep neural network surrogates and optimisation models in a unified mathematical framework. Based on this formulation recent research has increasingly explored the use of neural network surrogate models for the solution of stochastic programming problems by approximating computationally expensive recourse functions or simulation-based second stages. In particular, \cite{patel2022neur2sp} introduced Neur2SP, a framework where neural networks are trained to approximate the expected recourse function in two-stage stochastic programs, enabling near-optimal first-stage decisions with significantly reduced computation time. In \cite{Chou1} this approach has been proved to approximate well the recourse function of a stochastic facility location problem. 

More recently, \cite{Alcantara2025} proposed a Quantile Neural Network (QNN) approach to model the distributional characteristics of uncertain outcomes, offering improved performance in risk-averse formulations. Collectively, these studies show that neural network surrogates can substantially reduce computational costs while maintaining solution quality, positioning them as a promising direction for next-generation stochastic optimisation methods. Building on \cite{Chou1} we embed a deep neural network surrogate model to approximate the recourse function of an MHSP, aiming at finding solutions with enhanced performance on out-of-sample scenarios, showing the potential of surrogate models to generalise better on unseen scenarios.

\section{MHSP formulation}
\label{sec:MHSP formulation}
In this section, we establish the MHSP formulation. Typically, MHSP problems rely on continuous distributions to characterise both strategic and operational uncertainties. However, directly solving such infinite problems is often computationally intractable. Therefore, scenario tree approximations of the underlying stochastic processes are commonly employed to make the problem tractable. This is done by considering a finite number of realisations of both the random strategic processes and the operational processes.

The information structure at both the strategic and operational levels can be described in the form of a multi-horizon scenario tree, where at each strategic stage there is a discrete number of strategic nodes, where a specific realisation of the uncertain parameters at strategic level takes place. Let $\mathcal{I}$ be the set of ordered strategic nodes. To describe the operational uncertainty, we now consider operational subtrees in each strategic node $i \in \mathcal{I}$. Each strategic node $i$, except the root, is connected to a unique ancestor node at the previous stage, called its ancestor $a(i)$, and to nodes at the next stage, called its successors. 

Let $\Omega_i$ denote the set of possible operational scenarios at strategic node $i$. We indicate with $\omega_i$ the operational scenario derived by node $i \in \mathcal{I}$ and with $\pi^{\omega}_{i}$ its probability. 
Moreover, $\sum\limits_{\omega \in {\Omega}_i} \pi^\omega_{i}=1$, $ i \in \mathcal{I}$. Additionally, let $c_{i}$, $h_{i}$, $T_{i}$, $W_{i}$ be the uncertain vectors and matrices at strategic node $i \in\mathcal{I}$. If $i=1$, we assume $T_{1}=A$, $W_{1}=0$ (i.e., the null matrix), and $c_{1}$ and $h_{1}$ are known vectors.

Operational vectors and matrices at operational stage $t$, in operational scenario $\omega$ derived by node $i \in\mathcal{I}$, are given by $q_{i}^{\omega,t}$, $h_{i}^{\omega,t}$, $T_{i}^{\omega,t}$, and $W_{i}^{\omega,t}$. The strategic decision variable is given by $\mathbf{x}:=\left\{x_{i}\ |\ i \in\mathcal{I}\right\}$, with $x_{i}\in\mathbb{R}^{n_i}_{+}\times\mathbb{Z}_+^{n'_i}$. 
The operational decision variable is $\mathbf{y}:=\left\{y_{i}^{\omega,t}\ | \omega \in \Omega_i, t\in\mathcal{T}_i, i \in \mathcal{I}\right\}$, with 
$y_{i}^{\omega,t} \in\mathbb{R}^{n^{t}_{i}}_{+}\times\mathbb{Z}_+^{{n'}^t_i}.$ 
The MHSP is formulated as follows:
\begin{subequations}
\label{eq:MHSP}
\begin{alignat}{3}
\textit{MHSP}: \quad
    & \min_{\mathbf{x},\mathbf{y}}  
    && \sum_{i \in \mathcal{I}} \pi_{i} \Big( c_i^{\top} x_i 
       + \sum_{\omega \in \Omega_i} \pi^{\omega}_{i} 
         \sum_{t \in \mathcal{T}_i}  {q_{i}^{\omega,t}}^{\top} y_{i}^{\omega,t} \Big)
    \label{eq:1a}\\[4pt]
    & \text{s.t.}
    && Ax_{1} = h_{1},
    \label{eq:1b}\\[4pt]
    &
    && T_{i} x_{a(i)} + W_{i} x_{i} = h_{i},
    & i \in \mathcal{I}\setminus\{1\}, \label{eq:1c}\\[4pt]
    &
    && T_{i}^{\omega,1} x_{i} + W_{i}^{\omega,1} y_{i}^{\omega,1} = h_{i}^{\omega,1},
    & \omega \in \Omega_i,\; i \in \mathcal{I}, \label{eq:1d}\\[4pt]
    &
    && T_{i}^{\omega,t} y_{i}^{\omega,t-1} + W_{i}^{\omega,t} y_{i}^{\omega,t} = h_{i}^{\omega,t},
    & \omega \in \Omega_i,\; t \in \mathcal{T}_i\setminus\{1\},\; i \in \mathcal{I}. \label{eq:1e}
\end{alignat}
\end{subequations}

The objective function \eqref{eq:1a} minimises the expected total cost of strategic and operational decisions over the entire planning horizon.
Equation \eqref{eq:1b} defines a deterministic constraint on the first-stage strategic decisions. Equation \eqref{eq:1c} introduces a constraint that links strategic decisions across consecutive strategic stages. Constraint \eqref{eq:1d} establishes the connection between strategic decisions at node $i$ with operational decisions in the first operational period $t=1$.
Finally, constraint \eqref{eq:1e} connects operational decisions across consecutive operational periods within the same strategic stage.

\section{MHSP with surrogate neural network embedding}
\label{sec:surrogateNN_MHSP}
In this section, we mathematically present the embedding of a deep learning model as a surrogate for operational subproblems. We first introduce the deep learning networks in Section \ref{sec:deep learning networks} and we explain how to embed the deep learning surrogate into MHSP in Section \ref{sec: deep learning embedding}.

\subsection{Deep learning networks}
\label{sec:deep learning networks}
We use a deep learning model as a surrogate function for approximating the costs of the operational subproblems. In particular, we use a Feed Forward Neural Network (FFNN). An FFNN is a computational model that implements a parametric function 
\( f_{\theta}: \mathcal{X} \to \mathcal{Y} \), mapping an input \( x \in \mathcal{X} \) 
to an output \( y \in \mathcal{Y} \), where \( \theta \) denotes the set of learnable parameters. The network consists of multiple layers of neurons, also called units, each of which computes a parametric transformation of its inputs and produces an output that serves as input to the subsequent layer.

This transformation is typically non-linear and can be expressed as:
\begin{center}
    $h^{(l)} = \sigma \left( W^{(l)} h^{(l-1)} + b^{(l)} \right), \quad l = 1, \dots, L,$
\end{center}
where $h^{(0)} = x$ denotes the input vector, $h^{(l)}$ represents the activation of the $l$-th layer, and $W^{(l)}$ and $b^{(l)}$ are the learnable weights and biases of the layer, respectively. The function $\sigma$ is a non-linear activation function, and the output of the final layer, $h^{(L)}$, corresponds to the model's prediction $\hat{{y}}$.

A widely used choice for the activation function $\sigma$ is the Rectified Linear Unit (ReLU) (\cite{Nair2010RectifiedLU}), defined as the positive part of the affine transformation $a=W^{(l)} h^{(l-1)} + b^{(l)}$:
\begin{center}
    $\mathrm{ReLU}(a) = \max(0, a).$
\end{center}

Training a deep learning model involves optimising a suitable loss function  $\mathcal{L}( \hat{y}, y)$ with respect to the parameters  $\theta$ using gradient-based optimisation algorithms, most commonly Stochastic Gradient Descent (SGD) and its variants. 

A FFNN with a ReLU activation function can be embedded into a linear optimisation model by explicitly reformulating each neuron using auxiliary continuous and binary variables \citep{Fischetti2018}. For each neuron $h^{(l)}_j$ at position $j$ in layer $l$, its value can be represented as 
\begin{align}
h_j^{(l)} = \delta\left(\sum_{i=1}^{d^{(l-1)}} w^{(l-1)}_{ij} \, h^{(l-1)}_i + b^{(l-1)}_j \right),
\end{align}
where $d^{(l-1)}$ is the dimension of the previous layer, $w^{(l-1)}_{ij}$ and $b^{(l-1)}_j$ are the pretrained weights and bias, respectively, $\delta(\cdot)$ is the ReLU function.

To linearise the ReLU function, we introduce non-negative variables $\hat{h}_j^{(l)},\,\check{h}_j^{(l)} \ge 0$ and a binary activation indicator $z_j^{(l)} \in \{0,1\}$, and impose:
\begin{subequations}
\begin{align}
\sum_{i=1}^{d^{(l-1)}} w^{(l-1)}_{ij} \, \hat{h}^{l-1}_i + b^{l-1}_j&= \hat{h}_j^{(l)} - \check{h}_j^{(l)} \label{eq:neuron1}\\
z_j^{(l)} = 1 &\Rightarrow \hat{h}_j^{(l)} = 0 \label{eq:relu1}\\
z_j^{(l)} = 0 &\Rightarrow \check{h}_j^{(l)} = 0 \label{eq:relu2}
\end{align}
\end{subequations}

This set of constraints reproduces the piecewise linear behaviour of the ReLU activation: when the value on the left-hand side of equation~(\ref{eq:neuron1}) is strictly positive, the model forces $z_j^{(l)} = 0$ and $h_j^{(l)} = \hat{h}_j^{(l)} > 0$; when the left-hand side is non-positive, the model forces $z_j^{(l)} = 1$ and $\hat{h}_j^{(l)} = 0$. This approach ensures a complete representation of the neural network, and the number of constraints and variables introduced depends on the neural network size.  

\subsection{Deep Learning as a surrogate model} 
\label{sec: deep learning embedding}
We now extend the framework proposed in \cite{patel2022neur2sp} for two-stage stochastic programming models to account for the MHSP case.
In MHSP, the complicating variables are the strategic decisions $\mathbf{x}$ that link all the decision nodes. By fixing the complicating variables $\mathbf{x}$, the operational problems become independent. For a given strategic node $i$, the operational subproblem is formulated as 
\begin{subequations}
\label{eq:benders_subproblem_general}
\begin{alignat}{4}
    & \min_{ y_i^{\omega,t}}
    && \sum_{\omega \in \Omega_i} \pi^{\omega}_{i} 
       \sum_{t \in \mathcal{T}_i} {q_{i}^{\omega,t}}^{\top} y_{i}^{\omega,t}
    \\[4pt]
    & \text{s.t.}
    && T_{i}^{\omega,1} x_{i} + W_{i}^{\omega,1} y_{i}^{\omega,1} = h_{i}^{\omega,1},
    & \omega \in \Omega_i,\; i \in \mathcal{I},
    \\[4pt]
    &
    && T_{i}^{\omega,t} y_{i}^{\omega,t-1} + W_{i}^{\omega,t} y_{i}^{\omega,t} = h_{i}^{\omega,t},
    & \phantom{1111} \omega \in \Omega_i,\; t \in \mathcal{T}_i \setminus \{1\},\; i \in \mathcal{I}.
\end{alignat}
\end{subequations}

The single subproblem can include multiple operational scenarios. If the operational scenarios are independent of each other then an alternative approach would be to treat each as an independent subproblem. This has the potential advantage of maintaining a more accurate model of the operational problems and making it more efficient when creating training sets for neural networks. In this paper, we propose to use surrogate neural networks to approximate subproblem \eqref{eq:benders_subproblem_general}.

In particular, consider an operational scenario set $\Omega_i$. The expected value of the recourse function in problem (\ref{eq:1e}) can be expressed as:
\begin{equation}
\sum_{\omega \in \Omega_i} \pi^{\omega}_{i} \sum_{t \in \mathcal{T}_i} {q_{i}^{\omega,t}}^{\top} y_{i}^{\omega,t}.
\label{eq:recourse}
\end{equation}

This expected recourse function can be approximated by an FFNN by adapting the approach proposed in \cite{patel2022neur2sp}. \cite{patel2022neur2sp} considered two architectures for two-stage stochastic programming problems with integer first-stage decisions. The architecture consists of two networks: the first network generates an embedding of the scenario set, while the second network approximates the recourse function based on the concatenation of the first-stage solution and the embedded scenario representation.

In our setting, the first-stage decision variables are continuous. Consequently, constructing a training set that comprehensively covers the joint domain of first-stage solutions and possible scenarios would be computationally prohibitive. To address this challenge, we consider a fixed number of scenarios, which eliminates the need for scenario embeddings and reduces the dimensionality of the FFNN input. This approach enables the generation of a training set that sufficiently spans the feasible region of the first-stage problem while remaining computationally tractable. Formally, let the FFNN, $\hat{Q}_\theta({x}_i)$, approximate the expected value of the recourse function \eqref{eq:recourse}, where $\theta$ denotes the network parameters optimised during training. The input to the network is only the first-stage decision vector ${x}_i$ and the output is the predicted expected recourse value with respect to the scenario set $\Omega_i$. This methodology ensures that the neural approximation captures the dependency of the recourse function on first-stage decisions while maintaining a manageable network complexity.

We now describe the surrogate MIP for the learning model of problem \eqref{eq:1a}-\eqref{eq:1e} as follows:
\begin{subequations}
\label{eq:MHSP_nn}
\begin{alignat}{4}
    & \min_{\mathbf{x}}
    && \sum_{i \in \mathcal{I}} \pi_{i}\left(c_{i}^{\top} x_{i} + \hat{Q}_{\theta}(x_i)\right)
    \\[4pt]
    & \text{s.t.}
    && Ax_{1} = h_{1},
    \\[4pt]
    &
    && T_{i} x_{a(i)} + W_{i} x_{i} = h_{i},
    & i \in \mathcal{I}\setminus\{1\},
    \\[4pt]
    &
    && \sum_{k=1}^{d^{(l-1)}} 
       w^{(l-1)}_{kj}\,\hat{h}_k^{(l-1)} + b^{(l-1)}_j 
       = \hat{h}^{(l)}_{j} - \check{h}_{j}^{(l)},
    & \phantom{abc} l = 1,\ldots,L;\; j = 1,\ldots,d^{(l)},
    \label{eq:nn_1}
    \\[4pt]
    &
    && \sum_{k=1}^{d^{(L)}} w^{(L)}_{kj}\,\hat{h}_{j}^{(L)} + b^{(L)}_j 
       \leq \hat{Q}_{\theta}(x_i),
    & j = 1,\ldots,d^{(L)},
    \\[4pt]
    &
    && z^{(l)}_{j} = 1 \Rightarrow \hat{h}_{j}^{(l)} = 0,
    & l = 1,\ldots,L;\; j = 1,\ldots,d^{(l)},
    \\[4pt]
    &
    && z^{(l)}_{j} = 0 \Rightarrow \check{h}_{j}^{(l)} = 0,
    & l = 1,\ldots,L;\; j = 1,\ldots,d^{(l)},
    \\[4pt]
    &
    && \hat{h}^{(l)}_{j},\; \check{h}^{(l)}_{j} \ge 0,
    & l = 1,\ldots,L;\; j = 1,\ldots,d^{(l)},
    \\[4pt]
    &
    && z^{(l)}_{j} \in \{0,1\},
    & l = 1,\ldots,L;\; j = 1,\ldots,d^{(l)},
    \label{eq:nn_5}
\end{alignat}
\end{subequations}
where constraints \eqref{eq:nn_1}-\eqref{eq:nn_5} are aimed at emulating the ReLU activation signal $f(h) = \max\{0, Wh+b\}$ of each node $(i,j)$, $i=1\blue{,\ldots,}d^{(l-1)}$ and $j=1\blue{,\ldots,}d^{(l)}$, of the embedded NN model with $L$ layers and weight matrix $W$. Note that a binary variable $z^{(l)}_{j}$ must be introduced for each neuron of the neural network.

\section{Problem description, scenario generation and modelling assumptions}
\label{sec:problem_description_modelling_assumptions}
In this paper, we apply the proposed approach to solve a power system planning problem. In this section, we present the problem description and modelling strategies, including scenario generation and the modelling assumptions. 

The problem under consideration aims to choose (a) the optimal strategy for investment planning, and (b) operational scheduling for a power system to achieve emission targets at minimum overall costs under short-term uncertainty, including renewable energy availability and load profile, and long-term uncertainty, including CO$_2$ budget. 

For the investment planning, we consider: (a) thermal generators (Coal-fired plant, OCGT, CCGT, Diesel, nuclear plants); (b) generators with Carbon Capture and Storage (CCS) (Coal-fired plant with CCS and advanced CCS, gas-fired plant with CCS and advanced CCS); (c) renewable generators (onshore and offshore wind, solar); (d) electric storage (hydro pump storage and lithium); the capital expenditures and fixed operational costs coefficients are assumed to be known. 

The problem is to determine: (a) the capacities of technologies, and (b) operational strategies that include scheduling of generators, storage and approximate power flow among regions to meet the power demand with minimum overall investment and operational and environmental costs.

\subsection{Scenario generation and modelling assumptions}
In this section, we present the scenario generation and assumptions we use in the power system investment planning problem.

This paper adopts expert opinion methods for generating scenarios for long-term uncertain parameters. Given the long-term perspective, such data are typically derived from expected trends, historical observations, or expert projections. Consequently, we consider scenarios based on expert forecasts and projections. For short-term uncertainty, we adopt the random sampling approach from \cite{zhang2025multi-timescale}. Short-term scenarios are generated using the random sampling procedure proposed in \cite{zhang2025multi-timescale}. In the first step, the algorithm randomly selects an entire year from the dataset, which serves as the reference sample. This selected year is then divided into the four seasons. From each season, the algorithm extracts a sequence of equal length, sampling consecutive hours, preserving the temporal correlations and statistical properties of the original time series. Additionally, the algorithm incorporates one or more peak seasons, obtained by sampling consecutive time steps from periods of particularly high demand.

We assume that the Kirchhoff voltage law is omitted and we use a linear direct current power flow model.

\section{Power system model}
\label{sec:model}
This section presents the power system planning and operational optimisation model. The model is adapted from \cite{Zhang2024AUncertainty}.

\begin{description}[itemsep=-5pt, leftmargin=!,labelwidth=\widthof{\bfseries $\mathcal{I}^{Inv}_{j}xxxxx$}]

\item[\textbf{Investment planning model sets}]
\item[$p \in \mathcal{P}$] set of all technologies $(\mathcal{P} = \mathcal{G} \cup \mathcal{S} \cup \mathcal{R})$
\item[$i \in \mathcal{I}^{Ope}$] set of operational nodes
\item[$i \in \mathcal{I}^{Inv}$] set of investment nodes
\item[$i \in \mathcal{I}^{Inv}_{j}$] set of investment nodes $i$ $(i \in \mathcal{I}^{Inv})$ that are ancestors to operational node $j$ $(j \in \mathcal{I})$

\item[\textbf{Operational model sets}]
\item[$\omega \in \Omega$] set of operational scenarios
\item[$n \in \mathcal{N}$] set of time slices
\item[$t \in \mathcal{T}$] set of hours in all time slices
\item[$t \in \mathcal{T}_{n}$] set of hours belonging to time slice $n$
\item[$g \in \mathcal{G}$] set of thermal generators
\item[$s \in S$] set of electricity storage units
\item[$r \in \mathcal{R}$] set of renewable generators

\item[\textbf{Investment planning model parameters}]
\item[$C_{pi}^{Inv}$] unit investment cost of device $p$ in investment node $i$ $(p \in \mathcal{P}, i \in \mathcal{I}^{Inv})$ [£/MW]
\item[$C_{pi}^{Fix}$] unit fixed O\&M cost of device $p$ in operational node $i$ $(p \in \mathcal{P}, i \in \mathcal{I}^{Ope})$ [£/MW]
\item[$X_{pi}^{Hist}$] historical capacity of device $p$ $(p \in \mathcal{P}, i \in \mathcal{I}^{Ope})$ [MW]
\item[$X_{p}^{MaxAcc}$] maximum installed capacity of device $p$ $(p \in \mathcal{P})$ [MW]
\item[$X_{i}^{MaxInv}$] maximum built capacity in investment node $i$ $(i \in \mathcal{I}^{Inv})$ [MW]
\item[$\kappa$] scaling factor depending on operation years between investment nodes
\item[$\pi^{Inv}_{i}, \pi^{Ope}_{i}$] probability of investment node $i$ / operational node $i$ $(i \in \mathcal{I}^{Inv} \cup \mathcal{I}^{Ope})$
\item[$H_{p}^{P}$] lifetime of technology $p$ $(p \in \mathcal{P})$
\item[$X_{i}$] right-hand-side coefficients of operational subproblem $i$ $(i \in \mathcal{I}^{Ope})$
\item[$C_{i}$] cost coefficients of the operational subproblem $i$ $(i \in \mathcal{I}^{Ope})$
\item[$\mu_{i}^{E}$] $\mathrm{CO}_{2}$ budget at operational node $i$ $(i \in \mathcal{I}^{Ope})$
\item[$\mu_{i}^{D}$] scaling factor on power demand at operational node $i$ $(i \in \mathcal{I}^{Ope})$
\item[$S_{i}^{Ope}, S_{i}^{Inv}$] strategic stage of node $i$ $(i \in \mathcal{I}^{Ope} \cup \mathcal{I}^{Inv})$
\item[$C_{i}^{CO2}$] $\mathrm{CO}_{2}$ emission price at operational node $i$ $(i \in \mathcal{I}^{Ope})$

\item[\textbf{Operational model parameters}]
\item[$\pi^{\Omega}_{\omega}$] probability of operational scenario $\omega$ $(\omega \in \Omega)$
\item[$W_{t}$] scaling factor of operational period $t$ $(t \in \mathcal{T})$
\item[$H_{t}$] number of hours in operational period $t$ $(t \in \mathcal{T})$
\item[$\alpha_{g}^{G}$] maximum ramp rate of gas turbines $(g \in \mathcal{G})$ [MW/MW]
\item[$R_{r\omega t}^{R}$] renewable capacity factor for $r$ in scenario $\omega$, period $t$ $(r \in \mathcal{R}, \omega \in \Omega, t \in \mathcal{T})$
\item[$\eta_{s}^{S}$] storage efficiency $(s \in S)$
\item[$\gamma_{s}^{S}$] power ratio of storage unit $s$ $(s \in S)$ [MWh/MW]
\item[$E_{g}^{G}$] emission factor of generator $g$ $(g \in \mathcal{G})$ [tonne/MWh]
\item[$C_{g}^{G}, C_{s}^{S}$] operating unit cost of generator $g$ or storage $s$ $(g \in \mathcal{G}, s \in S)$ [£/MW]
\item[$C^{Shed}$] load shed penalty cost [£/MWh]
\item[$P_{\omega t}^{D}$] demand in scenario $\omega$, period $t$ $(\omega \in \Omega, t \in \mathcal{T})$ [MW]

\item[\textbf{Investment planning model variables}]
\item[$x_{pi}^{Acc}$] accumulated capacity of device $p$ in node $i$ $(p \in \mathcal{P}, i \in \mathcal{I})$ [MW]
\item[$x_{pi}^{Inv}$] newly invested capacity of device $p$ in investment node $i$ $(p \in \mathcal{P}, i \in \mathcal{I}^{Inv})$ [MW]

\item[\textbf{Operational model variables}]
\item[$p_{g}^{AccG}$] accumulated capacity of generator $g$ $(g \in \mathcal{G})$ [MW]
\item[$p_{r}^{AccR}$] accumulated capacity of renewable unit $r$ $(r \in \mathcal{R})$ [MW]
\item[$p_{s}^{AccS}$] accumulated charge/discharge capacity of storage $s$ $(s \in S)$ [MW]
\item[$p_{g\omega t}^{G}$] power generation of generator $g$ in scenario $\omega$, period $t$ $(g \in \mathcal{G}, \omega \in \Omega, t \in \mathcal{T})$ [MW]
\item[$p_{s\omega t}^{S+}, p_{s\omega t}^{S-}$] charge/discharge power of storage $s$ in scenario $\omega$, period $t$ $(s \in S,\omega \in \Omega,t \in \mathcal{T})$ [MW]
\item[$p_{\omega t}^{GShed}$] generation shed in scenario $\omega$, period $t$ $(\omega \in \Omega, t \in \mathcal{T})$ [MW]
\item[$q_{s\omega t}^{S}$] energy level of storage unit $s$ at start of period $t$ in scenario $\omega$ $(s \in S,\omega \in \Omega,t \in \mathcal{T})$ [MWh]
\item[$p_{\omega t}^{Shed}$] load shed in scenario $\omega$, period $t$ $(\omega \in \Omega,t \in \mathcal{T})$ [MW]

\end{description}

\subsection{Investment planning model}
The investment master problem Equations \eqref{eq:investment_cost}-\eqref{eq:mp_domain} follows the general formulation given by Equations \eqref{eq:1a}-\eqref{eq:1e}. The total cost for investment planning, Equation \eqref{eq:investment_cost}, consists of actual investment costs and the expected operational cost of the system over the time horizon $\kappa\sum_{i \in \mathcal{I}^{Ope}}\pi_{i}c^{Ope}_i$ which is the total approximated subproblem costs. Here, $\kappa$ is a scaling factor that depends on the time step between two successive investment nodes. 
\begin{equation}
    \min \sum_{i\in \mathcal{I}^{Inv}}\pi^{Inv}_{i}\sum_{p \in \mathcal{P}}C^{Inv}_{pi}x^{Inv}_{pi}+\kappa \sum_{i \in \mathcal{I}^{Ope}}\pi^{Ope}_{i}\sum_{p \in \mathcal{P}}C^{Fix}_{pi}x^{Acc}_{pi}+\kappa\sum_{i \in \mathcal{I}^{Ope}}\pi_{i}^{Ope}c^{Ope}_i ,\label{eq:investment_cost}
\end{equation}

Constraint \eqref{eq:cap_tech} states that the accumulated capacity of a technology $x^{Acc}_{pi}$ in an operational node equals the sum of the historical capacity $X^{Hist}_{p}$ and newly invested capacities $x^{Inv}_{pi}$ in its ancestor investment nodes $\mathcal{I}^{Inv}_{i}$ that are not retired.
\begin{equation}
    x^{Acc}_{pi}=X^{Hist}_{pi}+\sum_{j \in\mathcal{I}^{Inv}_{i}|\kappa(S^{Ope}_{i}-S^{Inv}_{j})\leq H^P_{p}}x_{pj}^{Inv}, \qquad p \in \mathcal{P},  i \in \mathcal{I}^{Ope}.\label{eq:cap_tech}
\end{equation}

Constraint \eqref{eq:max_built} ensures the maximum $X^{MaxInv}_{pi}$ capacity that is built in an investment node. Parameter $X^{MaxAcc}_{p}$ gives the maximum capacity that can be installed for different technologies.
\begin{equation}
    x^{Inv}_{pi} \leq X^{MaxInv}_{pi}, \qquad p \in \mathcal{P}, i \in \mathcal{I}^{Inv}. \label{eq:max_built}
\end{equation}

Constraint \eqref{eq:max_acc} establishes that the invested capacity and accumulated capacity of newly invested technologies and retrofitted technologies should be within the capacity limits. 
\begin{equation}
    x^{Acc}_{pi} \leq X^{MaxAcc}_{p},  \qquad  p \in \mathcal{P},  i \in \mathcal{I}^{Ope}, \label{eq:max_acc}
\end{equation}

The domains of variables are given as follows
\begin{equation}
   x^{Inv}_{pi}, x^{Acc}_{pi} \in \mathbb{R}^{+}_{0}. \label{eq:mp_domain}
\end{equation}

The vector $X_i=\left( \{x^{Acc}_{pi}, p \in \mathcal{P}\}, \mu^D_i, \mu^E_i \right)$ collects all right-hand side coefficients that will be fixed in operational subproblem, Equations \eqref{eq:SP_objective}-\eqref{eq:sp_domain}. The vector $C_i=\left(C^{CO_2}_{i}\right)$ collects all the cost coefficients. The vectors $X_i$ and $C_i$ will be fixed as parameters in the operational problem. 

\subsection{Operational problem}
We now present the operational problem in the deterministic equivalent model and note that we omit index $i$ and operational scenario $n$ in the operational model for ease of notation. 

The right-hand side parameters $P^{Acc}_p$, $P^{AccG}_g$, $P^{AccS}$, $Q^{AccS}_s$, $\mu^{D}$, and $\mu^E$ are fixed by the solution $X_i$ from solving the master problem Equations \eqref{eq:investment_cost}-\eqref{eq:mp_domain}. The CO$_2$ cost of generators that is included in parameter $C^G_g$ is fixed by $C_i$ from the master problem.

The operational cost $c^{Ope}(X_i,C_i)$ at one operational node $i$ is computed by solving subproblem Equations \eqref{eq:SP_objective}-\eqref{eq:sp_domain} given the decisions $x_i$ and $c_i$ made in Equations \eqref{eq:investment_cost}-\eqref{eq:mp_domain}. The objective function, the operational cost, includes total operating costs of generators $C^{G}_{g}p^{G}_{gt}$, power load shedding costs for power $C^{Shed}p_{t}^{Shed}$. $C^{G}_{g}$ includes the variable operational cost, fuel cost and the CO$_2$ tax, $C^{CO_2}$, charged on the emissions of generator $g$. 
\begin{equation}
    \min \sum_{t \in \mathcal{T}, \omega \in \Omega}\pi_{\omega}^{\Omega}W_t H_t\left(\sum_{g \in \mathcal{G}}C_g^{G}p_{g\omega t}^{S}+\sum_{s \in \mathcal{S}}C_s^{S}p_{s\omega t}^{G}+C^{Shed}p_{\omega t}^{Shed}\right). \label{eq:SP_objective}
\end{equation}

Constraints \eqref{eq:all_cap} ensure that the technologies operate within their capacity limits.
\begin{subequations}
\label{eq:all_cap}
    \begin{alignat}{3}
    &\quad&& p_{g\omega t}^{G}\leq P_{g}^{AccG}, & \qquad g \in \mathcal{G}, \omega \in \Omega, t \in \mathcal{T}, \label{eq:gen_cap}\\
    &\quad&& p_{s\omega t}^{S-} \leq P_{s}^{AccS}, & \qquad s \in \mathcal{S}, \omega \in \Omega, t \in \mathcal{T}, \label{eq:estore_discharge_cap}\\
    &\quad&& q_{s\omega t}^{S} \leq Q_{s}^{AccS}, & \qquad s \in \mathcal{S},  \omega \in \Omega, t \in \mathcal{T}, \label{eq:estore_energy_cap}
    \end{alignat}
\end{subequations}

Constraint \eqref{eq:ramp} captures how fast generators can ramp up or ramp down their power output, respectively. 
\begin{equation}
    -\alpha^{G}_{g}P_{g}^{AccG} \leq p_{g\omega t}^{G}-p_{g\omega, t-1}^{G}\leq \alpha^{G}_{g}P_{g}^{AccG}, \qquad  g \in \mathcal{G}, \omega \in \Omega, n \in \mathcal{N}, t \in \mathcal{T}_{n}. \label{eq:ramp}
\end{equation}

Constraint \eqref{eq:kcl} ensures that, in one operational period $t$, the sum of total power generation of generators $p^{G}_{g\omega t}$, power discharged from all the electricity storage $p^{SE-}_{s\omega t}$, renewable generation $R^{R}_{rwt}p^{AccR}_{r\omega t}$, power transmitted to this region, and load shed $p^{Shed}_{\omega t}$ equals the sum of power demand $\mu^{D}P^{D}_{\omega t}$, power transmitted to other regions, and power generation shed $p^{GShed}_{\omega t}$. The parameter $R^{GR}_{r\omega t}$ is the capacity factor of the renewable unit that is a fraction of the nameplate capacity $P^{AccR}$. 
\begin{equation}
\begin{split}
\sum_{g \in \mathcal{G}}p_{g\omega t}^{G}+\sum_{s \in \mathcal{S}}p_{s\omega t}^{S-}+\sum_{r \in \mathcal{R}}R^{R}_{r\omega t}P_{r}^{AccR}+p_{\omega t}^{Shed}=\mu^{D} P^{D}_{\omega t}+\sum_{s \in \mathcal{S}}p_{s\omega t}^{S+}+p_{\omega t}^{GShed}, \qquad \omega \in \Omega, t \in \mathcal{T}. \label{eq:kcl}
\end{split}
\end{equation}

Constraint \eqref{eq:storage_balance} states that the state of charge $q^{S}_{s\omega t}$ in period $t+1$ depends on the previous state of charge $q^{S}_{s\omega t}$, the charged power $\mu^{S}_sp^{S+}_{s\omega t}$ and discharged power $p^{S-}_{s\omega t}$. The parameter $\eta^{S}_{s}$ represents the charging efficiency. 

\begin{equation}
    q_{s\omega, t+1}^{S}=q_{s\omega t}^{S}+H_{t}(\eta_{s}^{S}p_{s\omega t}^{S+}-p_{s\omega t}^{S-}), \qquad  s \in \mathcal{S}, \omega \in \Omega, n \in \mathcal{N}, t \in \mathcal{T}_n. \label{eq:storage_balance}
\end{equation}

Constraint \eqref{eq:co2budget} restricts the total emission. The parameter $\mu^{E}$ is the CO$_2$ budget.
\begin{equation}
    \sum_{\omega \in \Omega}\pi_{\omega}^\Omega \sum_{t \in \mathcal{T}, g \in \mathcal{G}} W_tH_t E^{G}_{g} p^{G}_{g\omega t}\leq \mu^{E}. \label{eq:co2budget}
\end{equation}

The domains of variables are given as follows
\begin{equation}
\begin{split}
 p^{G}_{g\omega t}, p^{Shed}_{\omega t}, p^{GShed}_{\omega t}, p^{S+}_{s\omega t}, p^{S-}_{s\omega t}, q^{S}_{s\omega t} \in \mathbb{R}^{+}_{0}.\label{eq:sp_domain}
\end{split}
\end{equation}

\section{Case Study}
\label{sec:results}
In this section, we present the case study. We conduct a 15-year investment planning with a 5-year planning step. We implemented the algorithm and model in Python 3.11 using Pyomo and solved with Gurobi 12.0.2 \citep{gurobi}. The neural network is implemented with the TensorFlow framework version 2.19.0. The problem instances contain up to 3.43 million continuous variables and 6.62 million constraints. The training sets are created using a computer cluster with Dual socket AMD EPYC 9654 96-Core Processors which are equivalent to 192 cores and 650 GB RAM per node. The neural network training is done with a computer with Intel Core i7 - 1355 U with 10 cores and 1.7 GHz and 16 GB RAM. In addition, we choose to solve all the following instances to a 1\% convergence tolerance. 

\subsection{Data}
In this section, we introduce the data used in the computational study and describe how both long-term and short-term uncertainties are represented in the model. The data serve two purposes: first, to specify the parameters of the MHSP investment and operational models, and second, to generate training datasets for the neural network surrogate.

\subsubsection{Uncertainty data}
\begin{figure}[!ht]
\centering

\begin{subfigure}{0.48\textwidth}
    \centering
    \includegraphics[width=\textwidth]{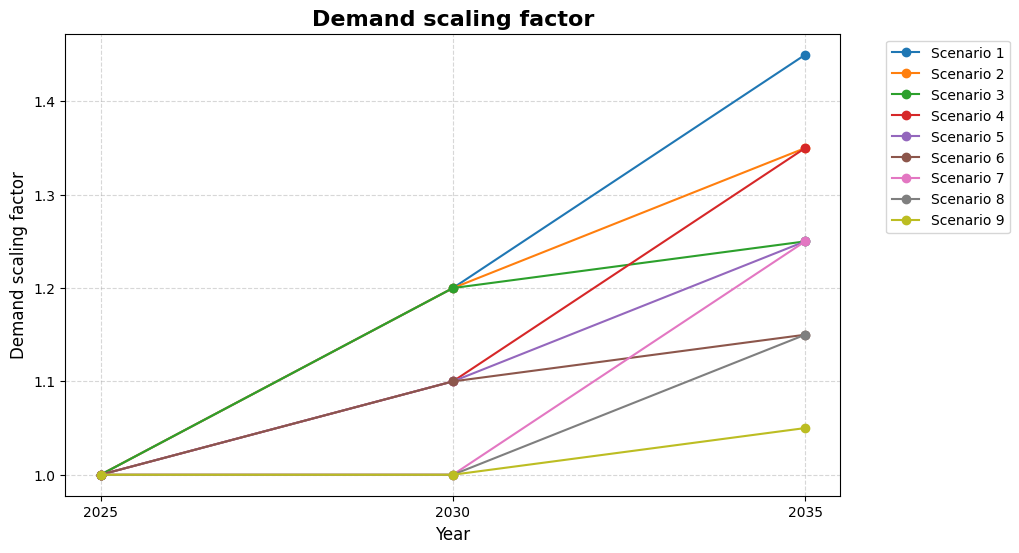}
    \caption{Long-term scenarios for demand scaling}
    \label{fig:4.4}
\end{subfigure}
\hfill
\begin{subfigure}{0.48\textwidth}
    \centering
    \includegraphics[width=\textwidth]{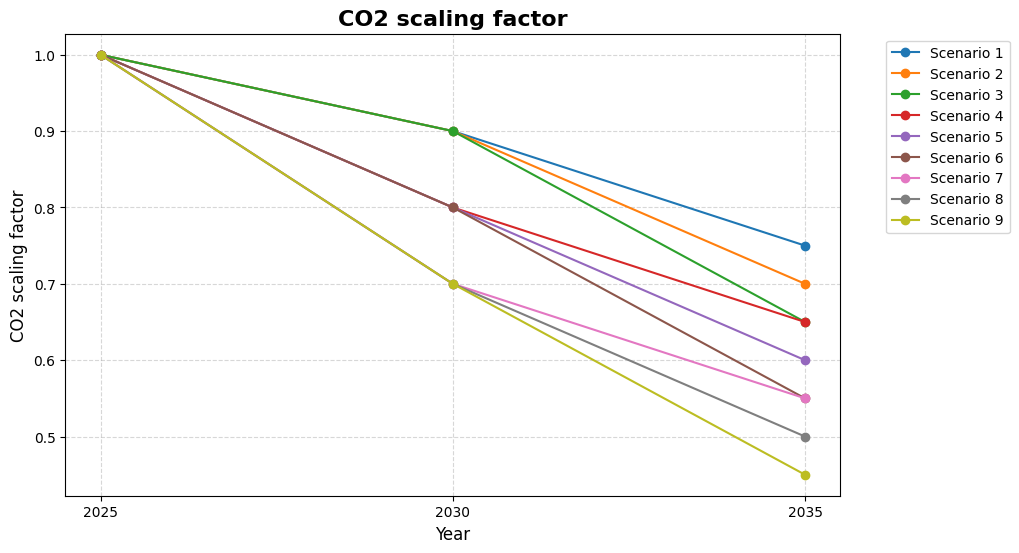}
    \caption{Long-term scenarios for CO$_2$ budget scaling}
    \label{fig:4.5}
\end{subfigure}

\caption{Long-term scenarios}
\end{figure}

In the case study, we consider long-term uncertainty in power demand scaling and the CO$_2$ limit scaling. A scaling factor is a multiplier that adjusts the baseline value of a parameter in different nodes. In practice, this means that the reference value is increased or decreased by applying the corresponding scaling factor at a specific node, in order to capture different possible evolutions of the parameter over time. The data used for the long-term uncertainty scenario tree is presented in Figures \ref{fig:4.4} and \ref{fig:4.5}. The long-term scenario tree includes 9 long-term scenarios and 13 nodes. The value of the uncertain parameters at each decision node is determined by the combination of the demand scaling and CO$_2$ scaling values corresponding to that node. 

We consider short-term uncertainty in time series parameters, including wind and solar capacity factors, and hourly power demand profile. We use the scenario generation method from \citep{zhang2025multi-timescale} to generate short-term scenarios. We have 5 years of data describing all the uncertain parameters, with a time frequency of 30 minutes. Using this algorithm, a set of short-term scenarios is generated using the same setting proposed in \cite{ZHANG2024106665}. This paper considers short-term scenarios composed of 96 time-step-length regular seasons and 25 time-step-length peak seasons, with one peak season in addition to the regular seasons. Using this setup, a set of 450 operational scenarios is generated to represent the underlying empirical distribution of our computational experiments.

\subsubsection{Creating training sets}
\label{sec:creating training sets}
To build the neural network surrogate model we need to generate a training set. The quality of a training dataset is an important factor for the performance of the neural network surrogate. A good training set should provide an adequate coverage of the data domain to ensure model generalisation capability, avoid overfitting and underfitting. The training dataset follows a typical regression structure: each sample consists of input features paired with a target value. In our case, the features are represented by the investment decisions $x^{Acc}_{pi}$ for every technology $p \in \mathcal{P}$, and the parameters $\mu^{D}_i$, $\mu^E_i$ and $C^{CO_2}_i$. The target is the operational cost in the respective decision node $i$. Consequently, every instance in the dataset is an artificially generated combination of features in a specific operational node, along with the associated operational cost.

To generate the input data for the neural network, a set of investment solutions is randomly generated. The random investment solutions generated should:
\begin{itemize}
    \item \textbf{Ensure data representativeness:} the dataset must capture the variability of the system and reflect a wide range of investment conditions, allowing the neural network to generalise effectively.
    \item \textbf{Achieve a good domain coverage:} since the decision variables are continuous, we must uniformly sample feasible points to obtain an adequate coverage of the data domain.
    \item \textbf{Guarantee solution feasibility:} each sampled point must satisfy the model’s constraints to ensure that the corresponding data are meaningful and consistent with real-world system behaviour.
\end{itemize}

For our problem we experimentally found that a Latin Hypercube Sampling (LHS) of the variables $x^{Inv}_{pi}$ with 50 intervals could satisfy the above conditions. To ensure feasibility $x^{Acc}_{pi}$ in each operational node were sampled within their feasible region, derived from (\ref{eq:cap_tech}) conditioned on the value of $x^{Inv}_{pi}$.

\subsection{Results}
In this section, we present the computational results from the case study. We begin by evaluating the neural network’s ability to approximate the recourse functions by measuring the approximation error across different neural network architectures. We then compare the solutions obtained from the surrogate model with those obtained by solving the deterministic equivalent, assessing both in-sample and out-of-sample stability for various sizes of the operational scenario set $\Omega$. Throughout the computational study, we use nine investment scenarios combined with five different sizes of $\Omega$, namely $|\Omega| \in \{5, 10, 15, 20, 50\}$. The complete information for the five cases is summarised in Table \ref{tab:case overview}, where monolithic model size refers to the number of variables and constraints contained in the deterministic equivalent model.

\begin{table}[!htb]
\caption{Overview of the cases used in the computational study}
\label{tab:case overview}
\centering
\resizebox{0.6\columnwidth}{!}{
\begin{tabular}{
    l
    S[table-format=3.0]
    S[table-format=1.2e1]
    S[table-format=1.2e1]
}
\toprule
 \textbf{Short-term}  & \textbf{Operational periods}  & \multicolumn{2}{c}{\textbf{Monolithic model size}} \\
 \textbf{scenarios}& \textbf{per short-term scenario} & \textbf{Continuous variables} & \textbf{Constraints} \\
\midrule
5   & 409  & 3.43e5 & 6.63e5 \\
10  & 409  & 6.87e5 & 1.32e6 \\
15  & 409  & 1.03e6 & 1.98e6 \\
20  & 409  & 1.37e6 & 2.65e6 \\
50  & 409  & 3.43e6 & 6.62e6 \\
\bottomrule
\end{tabular}
}
\end{table}

\subsubsection{Neural network approximation performance}
We now assess the ability of the neural networks to approximate the recourse costs of the operational subproblems. The recourse function (\ref{eq:SP_objective}) represents the expected operational cost, conditional on the investment decisions $x^{Inv}_{pi}, x^{Acc}_{pi} \in \mathbb{R}^{+}_{0}$ and with respect to the underlying set $\Omega$ of operational scenarios. For each investment node $i \in \mathcal{I}^{Inv}$, a neural network is trained to approximate the mapping from the vector of accumulated capacities and scaling parameters to the corresponding optimal operational cost.

Table \ref{tab:architectures} reports the three fully connected ReLU architectures, all with three hidden layers of increasing width: 16-8-4, 32-16-8, and 64-32-16 where the notation n$_1$-n$_2$-n$_3$ indicates the number of neurons in the first, second and third layer respectively. These architectures are trained and evaluated on problem instances of the operational model \eqref{eq:SP_objective}-\eqref{eq:sp_domain} with 5, 10, 15, 20 and 50 scenarios, respectively. The training data are generated by solving operational subproblems for investment decisions sampled via LHS, as described in Section \ref{sec:creating training sets}.

\begin{table}[!htb]
\centering
\caption{Neural network architectures.}
\resizebox{0.9\columnwidth}{!}{
\begin{tabular}{
    l
    S[table-format=4.0]
    S[table-format=4.0]
    S[table-format=4.0]
}
\toprule
\multicolumn{1}{c}{\textbf{Neural network structure}} &
\multicolumn{3}{c}{\textbf{Surrogate model size}} \\
\cmidrule(lr){2-4}
\textbf{Architecture} &
\textbf{Continuous variables} &
\textbf{Binary variables} &
\textbf{Constraints}
\\
\midrule
16-8-4   & 1643 &  336 & 1981 \\
32-16-8  & 2651 &  672 & 3325 \\
64-32-16 & 4667 & 1344 & 6013 \\
\bottomrule
\end{tabular}
}
\label{tab:architectures}
\end{table}

Table \ref{tab:nn_accuracy} summarises the approximation quality in terms of Mean Absolute Error (MAE), Mean Absolute Percentage Error (MAPE) and determination coefficient ($R^2$) for the three architectures and varying number of scenarios. Across all architectures and scenario set sizes, the $R^2$ values are consistently equal to 0.99, indicating that the surrogate model captures almost all of the variance in the operational cost. The MAPE values remain below 2.5\% in all tested configurations, and are typically close to 1\% for the largest network. As expected, larger architectures generally achieve lower MAE and MAPE, although the performance gains diminish beyond the 32-16-8 architecture, due to the well-known overfitting problems.

\begin{table}[!htb]
\centering
\caption{Results for neural network accuracy}
\label{tab:nn_accuracy}
\begin{tabular}{
    l
    S[table-format=4.2]
    S[table-format=1.2]
    S[table-format=1.2]
}
\toprule
\textbf{Architecture} & \textbf{MAE} & \textbf{MAPE (\%)} & \textbf{$R^2$} \\
\midrule

\multicolumn{4}{c}{\textbf{5 scenarios}} \\
\midrule
16-8-4   & 1290.39 & 2.40 & 0.99 \\
32-16-8  &  763.13 & 1.34 & 0.99 \\
64-32-16 &  641.65 & 1.11 & 0.99 \\[3pt]

\multicolumn{4}{c}{\textbf{10 scenarios}} \\
\midrule
16-8-4   &  977.27 & 2.17 & 0.99 \\
32-16-8  &  858.40 & 1.68 & 0.99 \\
64-32-16 &  612.01 & 1.05 & 0.99 \\[3pt]

\multicolumn{4}{c}{\textbf{15 scenarios}} \\
\midrule
16-8-4   & 1003.97 & 2.15 & 0.99 \\
32-16-8  &  818.55 & 1.64 & 0.99 \\
64-32-16 &  562.80 & 1.04 & 0.99 \\[3pt]

\multicolumn{4}{c}{\textbf{20 scenarios}} \\
\midrule
16-8-4   &  987.52 & 2.02 & 0.99 \\
32-16-8  &  837.77 & 1.63 & 0.99 \\
64-32-16 &  728.49 & 1.16 & 0.99 \\[3pt]

\multicolumn{4}{c}{\textbf{50 scenarios}} \\
\midrule
16-8-4   & 1249.90 & 2.28 & 0.99 \\
32-16-8  &  811.78 & 1.56 & 0.99 \\
64-32-16 &  669.64 & 1.00 & 0.99 \\
\bottomrule
\end{tabular}
\end{table}

These results confirm that a relatively compact FFNN is sufficient to approximate the recourse function with high accuracy, even when the number of operational scenarios increases. This supports the use of neural network surrogates as reliable replacements for the exact operational subproblems in the MHSP.

\subsubsection{Solution time}
We next compare the computational performance of the surrogate MHSP with that of the monolithic deterministic equivalent formulation. Table \ref{tab:solution time} reports, for each architecture and scenario size, the neural network training time, the MILP solving time for the surrogate model, their sum, and the solving time of the deterministic equivalent formulation. The last column reports the speed up factor, computed as the ratio between the deterministic solving time and the total time of the surrogate approach.

For small problem instances with 5 operational scenarios, the surrogate approach offers modest benefits: the 16-8-4 architecture yields a speed up of about 1.6, while the larger architectures are either comparable or slightly slower once training time is included. As the number of operational scenarios increases, however, the advantage of the surrogate-based approach becomes more significant. For 20 scenarios, the total solution time is reduced by factors between approximately 2.8 and 5, depending on the architecture. For the largest instance with 50 scenarios, the surrogate approach with the smallest network with a 16-8-4 architecture achieves a speed up of 34.72, while the deeper networks still provide an order of magnitude reduction in total computational time.

We also observe that, for the surrogate-based approach, the training phase typically dominates the total computational effort, whereas the solution of the surrogate-embedded MHSP is very fast. This suggests that, in settings where the same surrogate can be reused for multiple planning runs (e.g. under different policy scenarios or cost assumptions), the computational benefit of the surrogate approach would be even larger. Overall, these results indicate that neural network surrogates significantly improve scalability with respect to the number of operational scenarios.

\begin{table}[!htb]
\centering
\caption{Solving and training times of the surrogate models compared with the solution time of the deterministic equivalent (speed up: deterministic solving time divided by the total time of the surrogate models)}
\label{tab:solution time}
\resizebox{\textwidth}{!}{
\begin{tabular}{
    l
    S[table-format=3.2]
    S[table-format=3.2]
    S[table-format=3.2]
    S[table-format=4.2]
    S[table-format=1.2]
}
\toprule
\textbf{Architecture} &
\textbf{Training time (s)} &
\textbf{Solving time (s)} &
\textbf{Total (s)} &
\textbf{Deterministic solving time (s)} &
\textbf{Speed up}
\\
\midrule

\multicolumn{6}{c}{\textbf{5 Scenarios}} \\
\midrule
16-8-4   & 110.45 & 0.35  & 110.80 & 173.74 & 1.57 \\
32-16-8  & 130.92 & 7.64  & 138.56 & 173.74 & 1.25 \\
64-32-16 &  96.57 &102.49 & 199.06 & 173.74 & 0.87 \\
[4pt]

\multicolumn{6}{c}{\textbf{10 Scenarios}} \\
\midrule
16-8-4   & 204.37 & 0.40  & 204.77 & 484.32 & 2.37 \\
32-16-8  & 171.49 & 1.18  & 172.67 & 484.32 & 2.80 \\
64-32-16 & 129.50 &99.94  & 229.44 & 484.32 & 2.11 \\
[4pt]

\multicolumn{6}{c}{\textbf{15 Scenarios}} \\
\midrule
16-8-4   & 138.10 & 1.99  & 140.09 & 713.17 & 5.09 \\
32-16-8  & 223.02 &10.48  & 233.50 & 713.17 & 3.05 \\
64-32-16 & 215.43 &182.86 & 398.29 & 713.17 & 1.79 \\
[4pt]

\multicolumn{6}{c}{\textbf{20 Scenarios}} \\
\midrule
16-8-4   & 186.89 & 1.17  & 188.06 & 940.72 & 5.00 \\
32-16-8  & 208.90 & 4.88  & 213.78 & 940.72 & 4.40 \\
64-32-16 & 181.40 &153.46 & 334.86 & 940.72 & 2.81 \\
[4pt]

\multicolumn{6}{c}{\textbf{50 Scenarios}} \\
\midrule
16-8-4   &  84.89 & 1.78  &  86.67 & 3010.15 & 34.72 \\
32-16-8  & 267.14 & 6.60  & 273.74 & 3010.15 & 10.99 \\
64-32-16 & 154.97 &163.77 & 318.74 & 3010.15 &  9.44 \\
\bottomrule

\end{tabular}
}
\end{table}

\subsubsection{In-sample and out-of-sample stability}
\begin{figure}[htb!]
    \centering
    \includegraphics[scale=0.45]{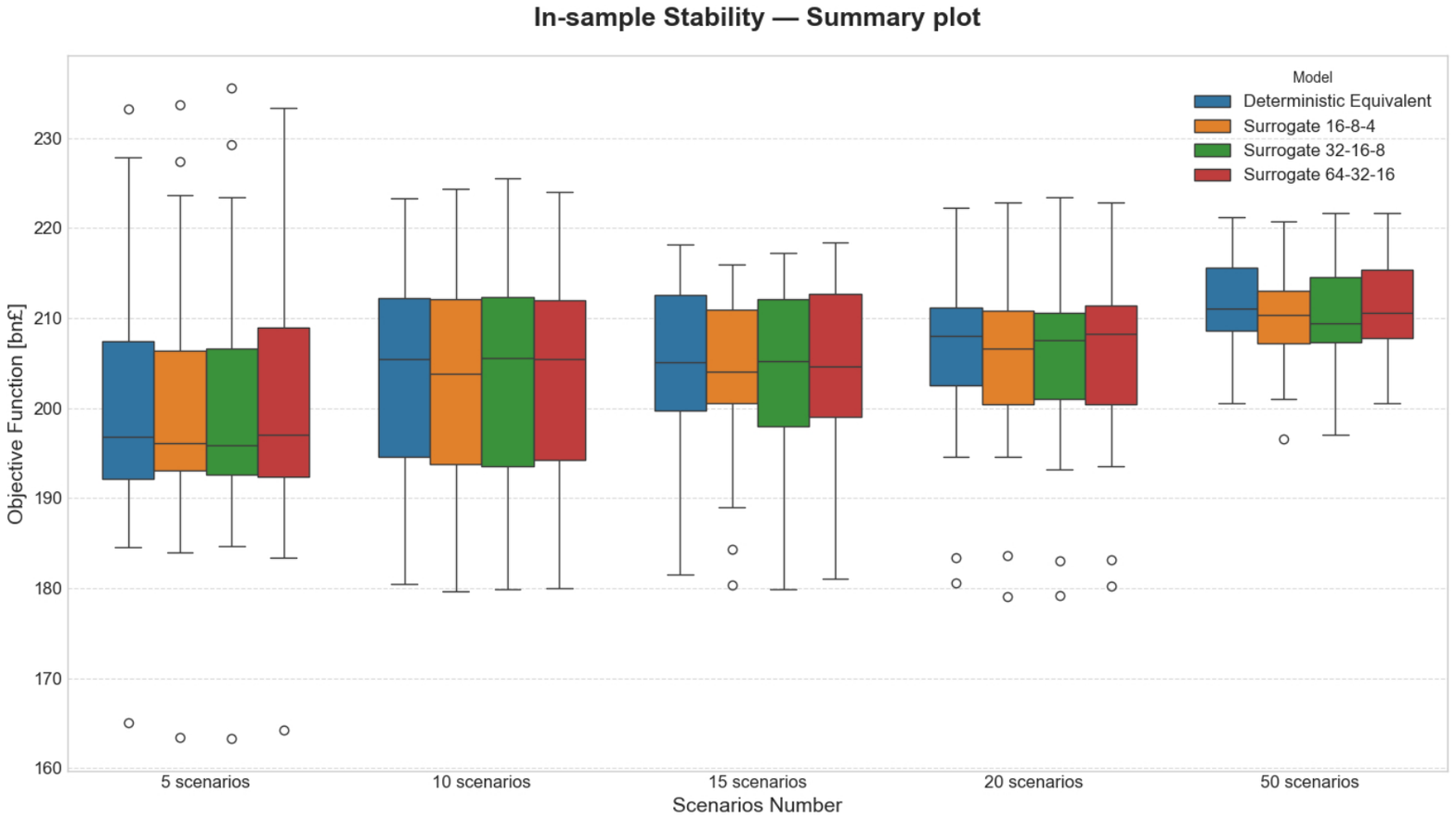}
    \caption{In-sample objective stability}
    \label{fig:is_stability}
\end{figure}

We finally investigate the robustness of the surrogate-based solutions by analysing their in-sample and out-of-sample stability. In-sample stability assesses how sensitive the optimal objective value is to the specific realisation of the scenario sample used in the optimisation model. For each problem instance, we generated 20 independent samples of operational scenarios from the same empirical distribution and solved both the deterministic equivalent and the surrogate-embedded MHSP. Figure~\ref{fig:is_stability} reports the distribution of the resulting objective values for varying scenario set sizes and architectures.

The box plots show that, across all scenario set sizes, the variability of the objective values is small for both approaches, indicating that the stochastic solutions are robust with respect to the sampled scenarios. The surrogate-based model exhibits a dispersion that is comparable to, and in several cases slightly smaller than, that of the deterministic equivalent formulation. This suggests that the approximation error introduced by the surrogate does not amplify the inherent sampling variability of the stochastic program.

To evaluate the out-of-sample stability, i.e. the generalisation properties of the solutions, we tested the solutions obtained from the deterministic equivalent model and the surrogate problem on an independent out-of-sample set of 200 operational scenarios. For each problem instance, the first-stage solutions from both the deterministic equivalent and the surrogate-embedded models were fixed, and their expected recourse costs evaluated over the out-of-sample scenarios.

\begin{figure}[htb!]
    \centering
    
    \begin{minipage}{0.4\textwidth} 
        \includegraphics[width=\textwidth]{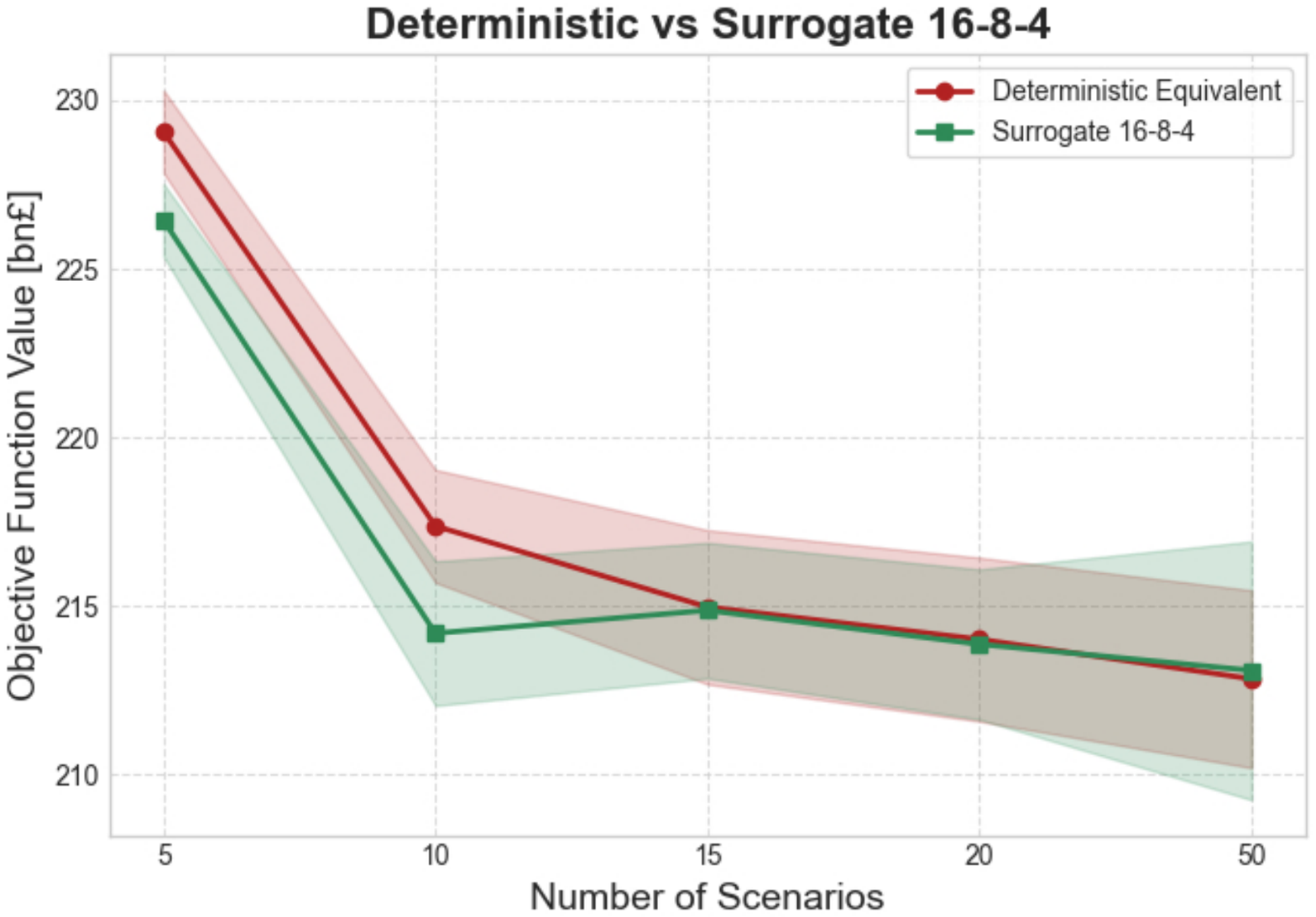}
        
    \end{minipage} 
    \begin{minipage}{0.4\textwidth}
        \includegraphics[width=\textwidth]{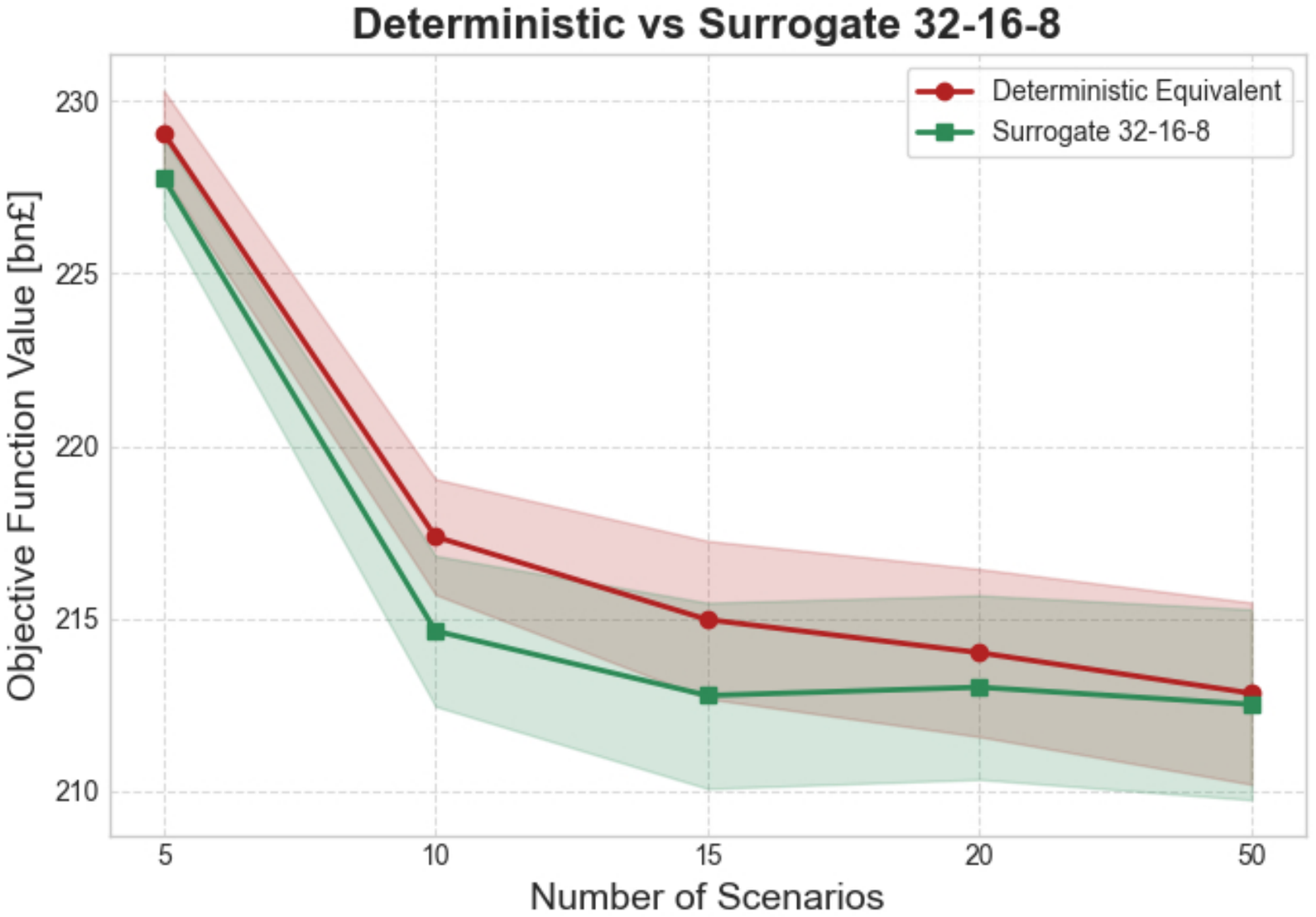}
        
    \end{minipage}
    
    \begin{minipage}{0.4\textwidth}
        \includegraphics[width=\textwidth]{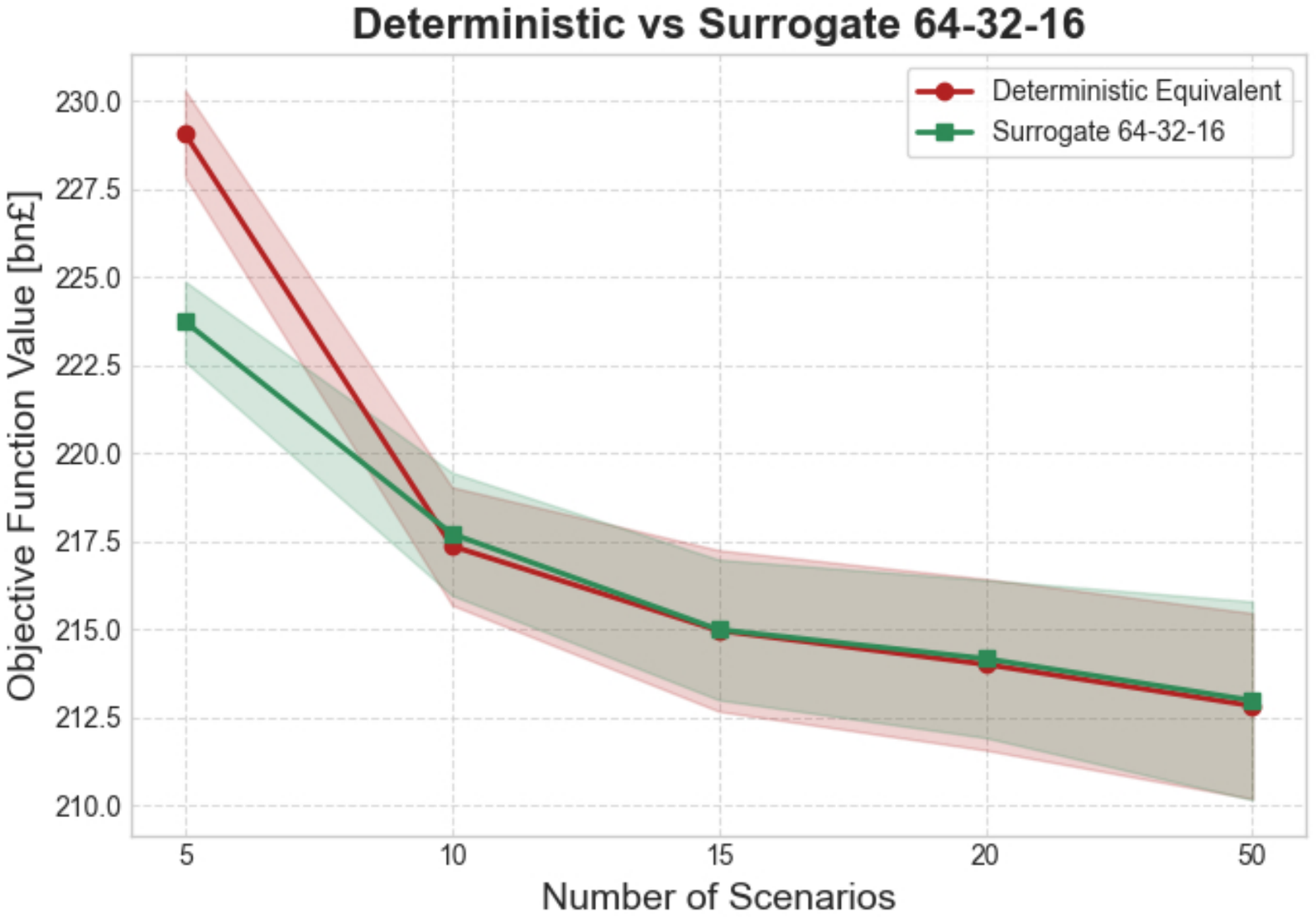}
        
    \end{minipage}
    \caption{Out-of-sample stability of the surrogate models and the deterministic equivalent}
    \label{fig:out-of-sample-stability}
\end{figure}

In Figure \ref{fig:out-of-sample-stability}, we report the mean, and the associated 95\% confidence intervals, of the optimal objective value (\ref{eq:investment_cost}) obtained by fixing the first-stage solutions derived from problems with different scenario set sizes. We observe that the generalisation ability of the surrogate model is comparable to that of the deterministic equivalent as the scenario set size increases. However, for small scenario sets, the first-stage solutions produced by the surrogate model consistently outperform those obtained from the deterministic equivalent, yielding lower cost values. This indicates that the surrogate neural network provides a better approximation of the recourse function and generalises more effectively when only a limited number of scenarios are available, highlighting the strength of machine learning approaches in capturing flexibility and enhancing generalisation.
Moreover, we notice that the problem where the embedded network architecture is 32-16-8 gives better results reaching a good trade-off between accuracy and generalisation ability. This confirms the hypothesis that the most effective surrogate model is not necessarily the one that most accurately approximates the recourse function value. This aligns with established principles in machine learning, where test accuracy alone is insufficient to guarantee generalisation capability. Indeed, overemphasis on accuracy can lead to overfitting, whereas models with better-calibrated uncertainties often support stronger generalisation and more stable downstream behaviour.

This observation highlights a broader methodological point: surrogate models should be evaluated not only on their pointwise approximation accuracy with respect to the corresponding recourse function value but also on their ability to preserve structural properties relevant to downstream tasks, in this case the optimisation problem. A surrogate model that slightly underperforms in estimating the exact value of the recourse function may nonetheless offer a representation that is more favourable to reach flexible solutions under conditions of uncertainty.
These insights suggest that surrogate-model selection should incorporate broader criteria such as robustness to perturbations, and alignment with the underlying task rather than relying solely on numerical approximation metrics. This shift in perspective can ultimately lead to more reliable recourse and improved flexibility of the decision pipeline as a whole. 

\section{Discussion}
\label{sec:discussion}
The computational study demonstrates that neural network surrogates can be effectively embedded into MHSP and used as accurate approximations of operational recourse costs. In this section we discuss the implications of these findings, the main limitations of the proposed framework, and possible directions for improvement.

From a methodological perspective, the results confirm that even relatively small FFNNs with ReLU activations are capable of capturing the dependency of the expected operational cost on the continuous first-stage decisions. This is non-trivial, as the recourse function arises from the optimal value of a large-scale linear programming or MILP which inherits its nonsmooth piecewise linear structure. The consistently high $R^2$ values and low percentage errors across all tested instances suggest that the surrogate preserves the essential shape of the recourse surface over the region of interest explored by the optimisation model. This explains why the surrogate-embedded MHSP recovers solutions with objective values and investment patterns that are very close to those of the deterministic equivalent model.

The main advantage of the proposed approach lies in its computational scalability. As the number of operational scenarios increases, the size of the deterministic equivalent formulation grows linearly, quickly leading to models with millions of variables and constraints. In contrast, the size of the surrogate-embedded model is essentially independent of the number of operational scenarios, as the scenarios only affect the offline training phase. The numerical results show that, for large scenario sets, this decoupling translates into substantial speed ups, while maintaining good solution quality and stability. This feature is particularly attractive in applications, such as energy system planning, where decision makers seek to explore many alternative assumptions and policy settings, and where richer scenario representations are needed to capture multi-timescale uncertainty. Moreover, the findings indicate that the surrogate neural network provides a reliable and computationally efficient alternative to the deterministic equivalent, offering particularly strong generalisation performance in settings with limited scenario availability.

The proposed framework has several limitations. Generating the training data requires solving a large number of operational subproblems, which can be computationally demanding. In the present study, this cost is acceptable and partially mitigated by parallel computation, but for more complex models or for very high-dimensional decision spaces this step could become a bottleneck. Techniques such as adaptive sampling, active learning, or the use of decomposition algorithms to accelerate subproblem solution may further reduce this burden.

Another limitation is that the embedded neural network relies on a piecewise linear MILP reformulation, which introduces additional continuous and binary variables. Although this was not a bottleneck in our experiments, the size of the surrogate representation grows with the number of layers and neurons, which may restrict the use of very large architectures. Alternative embedding strategies, such as cutting plane approximations, convex relaxations, or hybrid approaches that use surrogate models only in parts of the scenario tree, could offer an attractive trade-off between fidelity and size. Furthermore, we have assumed a fixed set of operational scenarios and did not investigate dynamically updating the surrogate when new data become available or when the scenario set is refined.

Finally, the case study focuses on a single country power system and uses a direct current power flow approximation, omitting some operational details such as network security constraints and unit commitment logic. While this level of detail is sufficient to illustrate the potential of the surrogate-based MHSP, further work is needed to assess the performance of the approach in more detailed operational models, in integrated multi-energy systems, and under alternative regulatory or market designs. Nevertheless, the results presented here provide encouraging evidence that neural network surrogates can be a valuable tool in MHSP, especially when computational tractability is a major concern.

\section{Conclusions and future work}
\label{sec:conclusions}
In this paper, we proposed a machine-learning-based framework for solving MHSP, and applied such a framework on power system investment planning under multi-timescale uncertainty. The key idea is to approximate the operational recourse function by training a feed-forward neural network on operational subproblem solutions, and to embed this surrogate directly into the MHSP formulation via a mixed-integer linear reformulation of the ReLU activations.

The case study on the UK power system shows that the surrogate model achieves high approximation accuracy across a range of neural network architectures and scenario set sizes. The surrogate-embedded MHSP produces objective values and structural characteristics that are very close to those obtained from the deterministic equivalent model. At the same time, the surrogate approach significantly reduces computational effort, especially when the number of operational scenarios is large, with up to 34.72 speed up. The in-sample and out-of-sample stability analyses further indicate that the surrogate-based solutions generalise well to unseen scenarios, and do not exhibit increased sensitivity to scenario sampling.

Future work includes (1) systematically comparing the proposed surrogate approach with decomposition algorithms for MHSP, such as Benders-type methods, in order to better understand the relative strengths and weaknesses of data-driven versus classical decomposition techniques, (2) exploring adaptive and active learning strategies for training data generation, in which the sampling of investment decisions is guided by the optimisation process itself, and (3) applying the proposed method to other problems with multi-timescale uncertainty, such as integrated energy systems, gas and hydrogen networks to further assess its practical impact and generality.

\section*{CRediT author statement}
\textbf{Hongyu Zhang:} Conceptualisation, Methodology, Validation, Formal analysis, Investigation, Visualisation, Data curation, Supervision, Writing - original draft, Writing - review \& editing. \textbf{Gabriele Sormani:} Conceptualisation, Methodology, Software, Validation, Formal analysis, Investigation, Visualisation, Data curation, Writing - original draft, Writing - review \& editing. \textbf{Enza Messina:} Conceptualisation, Methodology, Supervision, Formal Analysis, Writing - review \& editing, Funding acquisition. \textbf{Alan King:} Conceptualisation, Methodology, Supervision. \textbf{Francesca Maggioni:} Conceptualisation, Methodology, Validation, Formal analysis, Investigation, Visualisation, Writing - review \& editing. 

\section*{Declaration of competing interest}
The authors declare that they have no known competing financial interests or personal relationships that could have appeared to influence the work reported in this paper.





\clearpage

\setlength{\bibsep}{0pt plus 0.3ex}
\footnotesize{
\bibliographystyle{model5-names}
\bibliography{NN_MHSP}}
\clearpage

\appendix

\end{document}